\documentclass[10pt,twocolumn,twoside]{IEEEtran}
\usepackage{mathtools}
\usepackage{amsmath}
\usepackage{amssymb}
\usepackage{graphicx}
\usepackage{bbold}
\usepackage{enumerate}
\usepackage{enumitem}
\usepackage{framed}
\usepackage{color}
\usepackage{placeins}
\usepackage{tikz}
\usetikzlibrary{shapes, arrows, decorations, decorations.pathreplacing, decorations.pathmorphing, decorations.markings, fit, matrix}
\usepackage{algpseudocode}

\newtheorem{assumption}{Assumption}
\newtheorem{lemma}{Lemma}
\newtheorem{theorem}{Theorem}
\newtheorem{corollary}{Corollary}

\newtheorem{definition}{Definition}

\usepackage{soul}

\begin{document}
\title{LQG Control Performance with Low Bitrate Periodic Coding}
\author{       Behrooz~Amini,~\IEEEmembership{Student Member,~IEEE},  
        ~Robert~R~Bitmead,~\IEEEmembership{Fellow,~IEEE}
\thanks{The authors are  with the Department of Mechanical and Aerospace Engineering, University of California, San Diego, CA 92093-0411, USA, e-mail: \texttt{bamini,rbitmead@eng.ucsd.edu}.}}
\markboth{Submission to \textit{IEEE Transactions on Control of Network Systems}}
{Huang \MakeLowercase{\textit{et al.}}: Predictive Coding and Control}

\maketitle

\begin{abstract}
Specific low-bitrate coding strategies are examined through their effect on LQ control performance. By limiting the subject to these methods, we are able to identify principles underlying coding for control; a subject of significant recent interest but few tangible results. In particular, we consider coding the quantized output signal deploying period-two codes of differing delay-versus-accuracy tradeoff. The quantification of coding performance is via the LQ control cost. The feedback control system comprises the coder-decoder in the path between the output and the state estimator, which is followed by linear state-variable feedback, as is optimal in the memoryless case. The quantizer is treated as the functional composition of an infinitely-long linear staircase function and a saturation. This permits the analysis to subdivide into estimator computations, seemingly independent of the control performance criterion, and an escape time evaluation, which ties the control back into the choice of quantizer saturation bound. An example is studied which illustrates the role of the control objective in determining the efficacy of coding using these schemes. The results mesh well with those observed in signal coding. However, the introduction of a realization-based escape time is a novelty departing significantly from mean square computations.
\end{abstract}

\IEEEpeerreviewmaketitle

\section{Introduction}\label{sec:intro}
\IEEEPARstart{W}{e} consider a linear plant with input $u_t$ and output $y_t$ connected to a controller by a noise-free fixed-bitrate-$b$ memoryless channel. The measured output is coded for transmission through the channel and we consider several period-two coding or bitrate assignment strategies. In each case, the output is quantized with a linear fixed quantizer with saturation bound $\zeta$. The coding strategies perform a period-two bit-allocation for the signal being communicated across the channel. In Strategy~I, the $b$ bits of a $b$-bit quantizer are sent at each instant. Strategy~II applies a $2b$-bit quantizer and sends alternately the most significant $b$ bits and the least significant $b$ bits of the even-timed output sample only. The strategies differ in their delays and accuracy; $y_t$ has $b$ bits at each time versus $y_{2t}$ has $b$ bits at time $2t$ and $2b$ bits at time $2t+1$. No information is transmitted about $y_{2t+1}$ in the second strategy. A third, intermediate strategy is also examined. These coding/bit-assignment schemes are evaluated using the LQ performance of the controlled plant. Using a result from Curry \cite{Curry:1970}, the optimal control will comprise linear state-variable feedback and a conditional mean estimator using the decoded output.

\begin{figure}[h]
\begin{center}
\includegraphics[trim={270 11cm 0 11cm},clip,scale=.60]{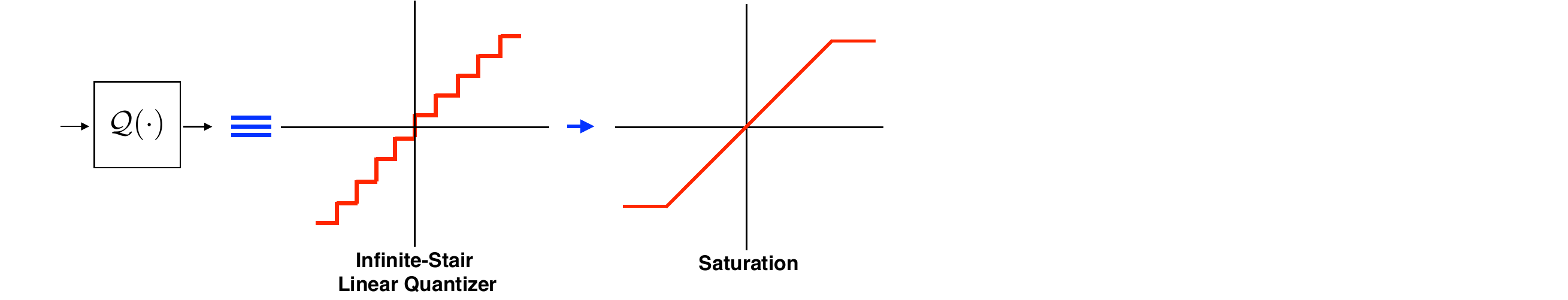}
\caption{Representation of a quantizer as the functional composition of two memoryless nonlinearities; an infinite quantizer and a saturation. The analysis treats each component in succession.\label{fig:Q2IS}}
\end{center}
\end{figure}

A quantizer is the functional composition or cascade of two distinct memoryless characteristic; an infinite quantizer and a saturation. This is depicted in Figure~\ref{fig:Q2IS}. 
We divide our analysis to consider each nonlinear aspect separately. In the case where the quantizer does not saturate and we use subtractive dithered quantization, the optimal conditional mean estimator is the Kalman filter, whose state estimate error covariance is computable using standard methods. The quantizer step size appears in the measurement noise term. For each coding strategy this covariance is simply computed and the LQ control performance derived. Next , the saturation nonlinearity is introduced by using these second-order signal statistics to compute the expected time before saturation. This \textit{escape time} is a function of the closed-loop controlled signal $y_t$. We use this property to tie to the controller objective function to the selection of the quantizer bound, $\zeta$. For a given feasible set of escape time, bitrate and control objective, there is a saturation bound and LQ performance. As the coding strategies change, so too does this performance. For a given escape time, we compare the quasi-stationary performance. 

\subsubsection*{Contribution of this paper}\mbox{\hskip 1mm}
\begin{itemize}
\item By treating a limited set of coding schemes, we are able to draw conclusions about coding in the output signal path. The range and the correlation/predictability of the closed-loop plant output play a role in the efficacy of coding. Less predictable outputs, such as those of minimum variance control, benefit less from coding. This concurs with observations in signal processing.
\item At low bitrates coding can become important.
\item The decomposition of the quantizer into two factors admits analysis using the escape time ideas introduced in this paper. This makes the study of methods possible by separating the estimator performance from the saturation behavior.
\item The escape time analysis permits the consideration of stabilization problems and performance together. The focus on realization based behavioral descriptors admits new viewpoints compared with asymptotic moments.
\end{itemize}

\subsubsection*{Relevant prior work}
Borkar and Mitter \cite{borkar1997lqg} study a full-state feedback formulation with vector quantization and coding delay similar to the strategies in this paper. They use the full state measurement to compute the process noise and then code this using vector quantization. They define a delay-accuracy tradeoff denoted by $N$ and indicating the number of noise samples held before transmission. Longer delay admits multiples of the underlying bitrate when eventually transmitted. When $N=1,$ their results are similar to our Strategy~I and when $N=2$ they resemble our Strategy~II. An optimal vector quantizer \cite{gray1984vector} then encodes the data into the available bits. This vector quantizer yields the conditional mean process noise reconstruction at the receiver even though the channel is not error-free. The decoded value is then used to construct the conditional mean state estimate. They could apply the same theorem from Curry in \cite{Curry:1970} to which we appeal shortly.  By limiting the discussion to stable systems, as in our earlier paper  \cite{huang2018predictive}, they are able to avoid explicit discussion of the saturation issues with quantization, vector or otherwise. 

Fu in \cite{fu2012lack} studied the coding problem for the control signal of linear quadratic Gaussian control with a memoryless error-free channel of fixed rate. The paper deals with optimization over the set of causal encoders and their decoder pairs. Fu looks only at finite-horizon optimal control and therefore sidesteps the stability and saturation questions. He does, however, develop a value for the finite-horizon LQG performance, which involves a distortion function that connects the coder and the objective function. In practice, optimizing this distortion appears intractable. He presents in Theorem 4.1, a corrected version of Fischer's result \cite{fischer1982optimal}, a \emph{weak} separation theorem where the estimate distortion function $D$ is a function of the control problem and not just the estimation problem parameters. When he considers the optimal coding problem for even a simple initial condition case, the solution depends on the control objective and the effect of the current encoding on future distortions. So his coder needs both memory and look-ahead and the problem begins to mimic the intractability of stochastic optimal control. In the current paper, by limiting our discussion to specific coding strategies, we reveal other aspects of a complicated picture. By restricting our coders to being memoryless, at least in Strategy I, we are able to appeal to the separation theorem of Curry.

Nair and Evans \cite{nair2004stabilizability} treat adaptive coding to achieve stabilization with limited bitrate. They assign on level of an adaptive quantizer to indicate saturation. When this level is received at the decoder, the quantizer range is expanded multiplicatively. Effectively, the bitrate required to stabilize an unstable system is tied to being able to achieve the expansion at a sufficiently rapid rate to catch the unstable output. This imaginative coding scheme concentrates on stabilization in mean square and does not address signal limits nor controlled performance.

The impact of quantization on performance at high rates is explored in \cite{gupta2006effect},  the state of the system being quantized prior to transmission to the controller, and they assess the performance of the controller 
to minimize a quadratic cost.
   
A similar approach is explored in \cite{fulton1997sampling} pertinent for speech coding but close to the current paper, in particular Strategies~I and II. These strategies are applied to speech with an autoregressive model. The performance is evaluated qualitatively by Mean Opinion Score. They show that down-sampling plus smoothing leads to better coding results for highly-correlated voiced speech and low-delay coding is preferable for unvoiced speech, which resembles modulated white noise.
      
Kostina and Hassibi \cite{kostina2019rate} consider LQ optimal control and the question of minimal channel capacity required to achieve a given bound on the expected LQ cost. They explore the problem with fully observed and partially observed state. In addition to this capacity bound, they explore specific lattice codes which achieve the bound. The bound itself depends on both control and estimation aspects for the partially observed case. They consider an error-free channel and explore all possible causal codes. Their communications structure is a limited capacity forward channel from the transmitter to receiver/controller together with a side channel which conveys the controller's state prediction back to the encoder. The minimizing codes transmit quantized versions of the error between the transmitter's state (or state prediction) and the receiver's state prediction. This communications structure obviates the requirement for the system to be stable. Although, similarly to \cite{nair2004stabilizability}, the logarithm of the determinant of the system matrix appears in the capacity bound.

\section{Problem statement}
Consider the following optimal control problem.
\begin{itemize}
\item Linear plant system with Gaussian noises:
\begin{align}
x_{t+1}&= Ax_t+Bu_t+w_t,\hskip 4mm x_0,\label{eq:nstate}\\
y_t&=Cx_t+v_t,\label{eq:noutput}
\end{align}
Here, state $x_t\in\mathbb R^n,$ input $u_t\in\mathbb R^p,$ output $y_t\in\mathbb R^m,$ process noise $w_t\in\mathbb R^n,$ measurement noise $v_t\in\mathbb R^m$. Noise sequences $\{w_t\}$ and $\{v_t\}$ are Gaussian, mutually independent, zero-mean and white with known covariances. The plant initial condition is also Gaussian and independent from $w_t$ and $v_t$ for all $t$.
\item Quadratic performance criterion, minimized over non-anticipatory controls, $u_t,$ computed from the received data at the controller,
\begin{align*}
J_{LQ}=\lim_{N \to \infty}\frac{1}{N}\text{E}\left(\sum_{t=1}^Nx_t^TQ_cx_t+u_t^TR_cu_t\right).
\end{align*}
\item The communications link between the plant measurement and the control computation consists of a limited bitrate, $b$-bits-per-sample, memoryless noise-free channel.
\item The coder-controller is restricted to the following elements. 
\begin{itemize}
\item The measurement $y_t$ is quantized to a fixed number of bits, which can be larger than $b$.
\item Some of these bits are encoded into the bitstream forwarded to the controller subject to the bitrate limit.
\item For this paper, we restrict attention to period-two bitrate assignment strategies.
\item The controller computes and applies the control.
\end{itemize}
\end{itemize}

\subsection{LQ Optimal Controller}
Denote by $\{p_t\}$ the sequence of decoded signal values available at the controller. Then, we have the following result from Curry.
\begin{theorem}[Curry \cite{Curry:1970}]\label{thm:curry}
For the linear state system (\ref{eq:nstate})-(\ref{eq:noutput}) with nonlinear, memoryless measurement 
\begin{align*}
p_t&=\varphi_t(x_t,v_t),
\end{align*}
with $\{v_t\}$ white and independent from $x_t$ and quadratic objective function
\begin{align}\label{eq:quadobj}
J_t=\text{E}\left(\left.\sum_{k=t}^{N+1}{x_k^TQ_kx_k+u_k^TR_ku_k}\right|\mathbf P^t,\mathbf U^{t-1},\pi_{0|-1}\right),
\end{align}
the optimal output feedback control is given by
\begin{align*}
u^\star_t=-K_t\,\text{E}(x_t|\mathbf P^t,\mathbf U^{t-1},\pi_{0|-1}),
\end{align*}
where, $\mathbf P^t=\{p_1,p_2,...,p_t\}$, $\mathbf U^t=\{u_1,u_2,...,u_t\}$, $\pi_{0|-1}$ is the initial state density, and $K_t$ is the  LQ optimal feedback gain.
\end{theorem}

Decoding signal $p_t$ at the receiver side, the filtered plant state estimate $\hat x_{t|t}$ and infinite-horizon control law  $u_t=-K\hat x_{t|t}$  are computed with $K=\text{dare}(A,B,Q_c,R_c).$ The performance is evaluated with the LQ  criterion. Signal $p_t$ will be derived from output $y_t$ by quantization and coding.

\section{Controller Coding Strategies}
The problem statement imposes the quantization of plant output signal $y_t$. We restrict our attention to uniform quantization and limit consideration to subtractive dithered quantizers in order to facilitate the receiver-side estimation.

\subsection{Dithered Quantization }
A \textit{subtractive b-bit dithered quantizer,} $\mathbf Q_b(\cdot),$ is a memoryless function which takes input signal $y_t$ and dither signal, $d_t,$ and produces an output signal
\begin{align}\label{eq:qdef}
\mathbf Q_b(y_t)=Q_b(y_t+d_t)-d_t,
\end{align}
where $Q_b(\cdot)$ is a standard uniform quantizer. Such quantizers are examined in detail in, for example, \cite{Widrow&Kollar08}.

\begin{theorem}\label{th:qtsd}
Consider a uniform, midrise, symmetric, $b$-bits-per-channel, subtractive dithered quantizer, $\mathbf Q_b(\cdot)$, with saturation bounds $\pm\zeta$. Assume:
\begin{enumerate}[label=(\Alph*)]
\item dither $d_t$ is a white noise process independent from $y_t$ with a probability density possessing characteristic function, $\Phi_d(\cdot),$ satisfying $\Phi_d\left(l\frac{\pi 2^b}{\zeta}\right)=0$ for $l=\pm 1,\pm 2, \dots$,
\item\label{ass:saturate} $y_t+d_t\in[-\zeta,\zeta],$ i.e. no saturation of the dithered quantizer occurs.
\end{enumerate}
Then, the quantization error
\begin{align}
q_{b,t}\triangleq \mathbf Q_b(y_t)-y_t,
\end{align}
is: (i) white, (ii) independent from $y_t$, (iii) uniformly distributed on $\left[-\frac{\zeta}{2^b},\frac{\zeta}{2^b}\right].$ 
\end{theorem}

This theorem, an embellishment of Theorem~QTSD of \cite{Widrow&Kollar08}, presents conditions under which the quantization error is an additive white noise independent from the signal being quantized as studied with details in \cite{gray1993dithered}. Denote the quantizer step size as 
\begin{align*}
\Delta=\frac{\zeta}{2^{b-1}}.
\end{align*}
Then, we note the following for this dithered quantizer.
\begin{align}
&q_{b,t}\sim\mathcal U\left[-\frac{\Delta}{2},\frac{\Delta}{2}\right],\;\;\text{E}(q_{b,t})=0,\nonumber\\
&\text{cov}(q_{b,t})=\frac{\zeta^2}{3\times 2^{2b}}\triangleq S_b.\label{eq:sbdef}
\end{align}
We also, note that the characteristic function condition is satisfied by dither which is uniformly distributed $\mathcal U[-\Delta/2,\Delta/2]$ or which is triangularly distributed $\text{tr}[-\Delta,\Delta]$, for example. In our calculations later, we use uniform dither $d_t$.

\subsection{Period-two Bit-assignment and Transmission Strategies}
We consider a fixed-rate, $b$-bits-per-transmission, channel and propose three period-two quantization strategies which reflect similar approaches from Signal Processing  \cite{goodwin2014use}, \cite{cea2014control}, \cite{fulton1997sampling}. The intention is to manage the quantization error with periodic changes to the effective bitrate and allied signal delay. We will examine the efficacy of these methods in terms of their benefits for LQ output feedback control.

The presence of the $b$-bits-per-sample channel militates that the subtractive dithered quantizer operates on both sides of the channel. That is, $b$ bits are transmitted each sample as symbol $m_t$ from the transmitter. Then at the receiver subtractive dither is applied. This and other implementation issues of wordlength etc. are discussed in \cite{Widrow&Kollar08}. With our period-two strategies, both the dithering and the subtraction will be modified. Here $\text{MSB}_n(x_t)$ and $\text{LSB}_n(x_t)$ denote the most significant and least significant $n$ bits of signal $x_t$. While $m_t$ is the $b$-bit transmitted message at time $t,$ $p_t$ or $p^\prime_t$ denotes the reconstructed/decoded plant output at the receiver for input into the Kalman filter.

\textit{Strategy I}\label{strat:one}
\begin{algorithmic}[1]
\For{$t$ even or odd}
\State $m_t=Q_{b}(y_t+d^b_t)$ is transmitted
\State $p_t\gets m_t-d^b_t$ at the receiver
\State $\hat x_{t|t}\gets$ \eqref{eq:sIall} Kalman filter Lemma~\ref{thm:anderson}
\State $u_t\gets -K\hat x_{t|t}$
\EndFor
\end{algorithmic}

\textit{Strategy II}
\begin{algorithmic}[1]
\If{$t=2k$, even time,}
\State $m_{2k}=\text{MSB}_b(Q_{2b}(y_{2k}+d^{2b}_{2k}))$ is transmitted
\State $p_{2k}\gets m_{2k}$ at the receiver without dither subtraction
\State $\hat x_{2k|2k}\gets$ \eqref{eq:sIIeven} Kalman filter from Lemma~\ref{lem:lem2}
\State $u_{2k}\gets -K\hat x_{2k|2k}$
\Else \;$t=2k+1$, odd time,
\State $m_{2k+1}=\text{LSB}_b(Q_{2b}(y_{2k}+d^{2b}_{2k}))$ is transmitted
\State $p_{2k+1}\gets p_{2k}+2^{-b}m_{2k+1}-d^{2b}_{2k}$ at the receiver
\State $\hat x_{2k+1|2k+1}\gets$ \eqref{eq:sIIodd} Kalman filter from Lemma~\ref{lem:lem2}
\State $u_{2k+1}\gets -K\hat x_{2k+1|2k+1}$
\EndIf
\end{algorithmic}

\textit{Strategy III}
\begin{algorithmic}[1]
\If{$t=2k$, even time,}
\State $m_{2k}=\text{MSB}_b(Q_{b+r}(y_{2k}+d^{b+r}_{2k}))$ is transmitted
\State $p_{2k}\gets m_{2k}$ at the receiver without dither subtraction
\State $\hat x_{2k|2k}\gets$ \eqref{eq:sIIIeven} Kalman filter from Lemma~\ref{lem:lem3}
\State $u_{2k}\gets -K\hat x_{2k|2k}$
\Else \;$t=2k+1$, odd time,
\State $m_{2k+1}=\text{LSB}_r(Q_{b+r}(y_{2k}+d^{b+r}_{2k}))$
\State\hskip 20mm$+2^{-r}\text{MSB}_{b-r}(Q_{b-r}(y_{2k+1}+d^{b-r}_{2k+1}))$ 
\State\hskip 10mm is transmitted
\State $p^\prime_{2k}\gets p_{2k}+2^{-b}\text{MSB}_r(m_{2k+1})-d^{b+r}_{2k}$
\State $p_{2k+1}\gets \text{LSB}_{b-r}(m_{2k+1})-d^{b-r}_{2k+1}$
\State $\hat x_{2k+1|2k+1}\gets$ \eqref{eq:sIIIodd} Kalman filter from Lemma~\ref{lem:lem3}
\State $u_{2k+1}\gets -K\hat x_{2k+1|2k+1}$
\EndIf
\end{algorithmic}
%
%
%

We note two central features of the time-varying strategies.
\begin{itemize}
\item[$\bullet$] Strategy~I uses a quantizer of step size $\frac{\zeta}{2^{b-1}}$ and associated dither $d^b_t\sim\mathcal U[-\zeta/2^b,\zeta/2^b]$,   while Strategy~II uses step size $\frac{\zeta}{2^{2b-1}}$ and $d^{2b}_{2k}\sim\mathcal U[-\zeta/2^{2b},\zeta/2^{2b}]$, and Strategy~III uses alternately $\frac{\zeta}{2^{2(b+r)-1}}$ and $\frac{\zeta}{2^{2r-1}}$ for the step size and the associated dithers. 
\item[$\bullet$] Strategies~II and III at even times receive undithered $b$-most-significant-bit transmissions, since the dither operates further along the bitstream. Accordingly, the quantization error at even times is not white, nor uniform, nor independent from $y_{2k},$ even though the quantization noise for $y_{2k}$ at time $2k+1$ does possess these properties. We shall conduct our analysis blithely without taking these even quantization error properties fully into account.
\end{itemize}


We note that, with Strategies~II and III, the state estimate calculation will be non-standard at the controller, reflecting the periodic information pattern. The associated Kalman filter will be presented shortly and computes $\hat x_{2k|2k+1}$ and then $\hat x_{2k+1|2k+1}$ from the received data. The derivation of these filters and their properties is a core contribution of the paper.

\section{Kalman Filters and Covariances for the Strategies}\label{sec:KF}
We derive the Kalman filters associated with each of the strategies under the following assumption.
\begin{assumption}[For this and the following sections alone]\label{ass:nosat}
The quantizer never saturates. That is, $y_t+d_t\in[-\zeta,\zeta]$. So, following Theorem~\ref{th:qtsd}, the quantization errors:
\begin{description}
\item[\rm Strategy I:] $p_t-y_t$;
\item[\rm Strategy II:] $p_{2k+1}-y_{2k}$;
\item[\rm Strategy III:] $p^\prime_{2k}-y_{2k}$ and $p_{2k+1}-y_{2k+1}$;
\end{description}
are independent from $\{y_t\}$, white, zero-mean, uniformly distributed with covariances $S_b$, $S_{2b}$, $S_{b+r}$ and $S_{b-r}$ respectively, where $S_b$ is defined in \eqref{eq:sbdef}.

Further and without justification, we assume that the other quantization errors, $p_{2k}-y_{2k}$ in Strategies II and III, satisfy the same properties with covariances $S_b$.
\end{assumption}
\begin{assumption}\label{ass:initial}
Each strategy's state estimator commences with state estimate $\hat x_{0|-1}$ and covariance $\Sigma_{0|-1}, \text{at } t=0.$
\end{assumption}
The following results for Strategies I, II and III are proved in the Appendix.
\begin{lemma}[Anderson, Moore \cite{anderson2012optimal}]\label{thm:anderson}
For Strategy I, the  Kalman filter driven by signal $p_t$ from Algorithm~I Line~3 is calculated by: 
\begin{align}
&L_{t}=\Sigma_{t|t-1}C^T(C\Sigma_{t|t-1}C^T+R+S_b)^{-1},\nonumber\\
&\hat x_{t|t}=\hat x_{t|t-1}+L_{t}(p_{t}-C\hat x_{t|t-1}),\label{eq:sIall}\\
&\hat x_{t+1|t}=(A-BK)\hat x_{t|t},\nonumber\\
&\Sigma_{t+1|t}=A\Sigma_{t|t-1}A^T-AL_tC\Sigma_{t|t-1}A^T+Q.\nonumber
\end{align}
\end{lemma}
\begin{lemma}\label{lem:lem2}
For Strategy II, the  Kalman filter driven by signals $p_{2k}$ from Algorithm~II Line~3 and $p_{2k+1}$ from Line~8 is calculated by:
\vskip 1mm
\noindent\textit{At even times, $t=2k:$}
\begin{align}
&L_{2k}=\Sigma_{2k|2k-1}C^T(C\Sigma_{2k|2k-1}C^T+R+S_b)^{-1},\nonumber\\
&\hat x_{2k|2k}= \hat x_{2k|2k-1}+L_{2k}(p_{2k}-C\hat x_{2k|2k-1}),\label{eq:sIIeven}
\end{align}
\vskip 1mm
\noindent\textit{At odd times, $t=2k+1:$}
\begin{align}
&L_{2k+1}=\Sigma_{2k|2k-1}C^T(C\Sigma_{2k|2k-1}C^T+R+S_{2b})^{-1},\nonumber\\
&\hat x_{2k+1|2k+1}= A\left(\hat x_{2k|2k-1}+L_{2k+1}(p_{2k+1}-C\hat x_{2k|2k-1})\right)\nonumber\\
&\hskip 58mm-BK\hat x_{2k|2k},\label{eq:sIIodd}\\
&\hat x_{2k+2|2k+1}=(A-BK)\hat x_{2k+1|2k+1},\nonumber\\
&\Sigma_{2k+2|2k+1}=A^2\Sigma_{2k|2k-1}{A^2}^T-A^2L_{2k+1}C\Sigma_{2k|2k-1}{A^2}^T\nonumber\\
&\hskip 60mm+AQA^T+Q.\nonumber
\end{align}
\end{lemma}
%
%

\begin{lemma}\label{lem:lem3}
For Strategy~III, the  Kalman filter driven by signals $p_{2k}$ from Algorithm~III Line~3, $p^\prime_{2k}$ Line~10 and $p_{2k+1}$ Line~11 is calculated by:
\vskip 1mm\noindent\textit{At even times, $t=2k,$}
\begin{align}
&L_{2k}=\Sigma_{2k|2k-1}C^T(C\Sigma_{2k|2k-1}C^T+R+S_b)^{-1},\nonumber\\
&\hat x_{2k|2k}= \hat x_{2k|2k-1}+L_{2k}(p_{2k}-C\hat x_{2k|2k-1}),\label{eq:sIIIeven}
\end{align}
\vskip 1mm\noindent
\textit{At odd times, $t=2k+1,$}
\begin{align}
&L^\prime_{2k}=\Sigma_{2k|2k-1}C^T(C\Sigma_{2k|2k-1}C^T+R+S_{b+r})^{-1},\nonumber\\
&\hat x^\prime_{2k+1|2k}=A\hat x_{2k|2k-1}+AL_{2k}^\prime(p^\prime_{2k}-C\hat x_{2k|2k-1})\nonumber\\
&\hskip 50mm-BK\hat x_{2k|2k},\nonumber\\
&\Sigma^\prime_{2k+1|2k}=A\Sigma_{2k|2k-1}A^T-A\Sigma_{2k|2k-1}C^T\times\nonumber\\
&\hskip 10mm(C\Sigma_{2k|2k-1}C^T+R+S_{b+r})^{-1}C\Sigma_{2k|2k-1}A^T+Q,\nonumber\\
&L_{2k+1}=\Sigma^\prime_{2k+1|2k}C^T(C\Sigma^\prime_{2k+1|2k}C^T+R+S_{b-r})^{-1},\nonumber\\
&\hat x_{2k+1|2k+1}=\hat x^\prime_{2k+1|2k}+L_{2k+1}(p_{2k+1}-C\hat x^\prime_{2k+1|2k}),\label{eq:sIIIodd}\\
&\hat x_{2k+2|2k+1}=(A-BK)\hat x_{2k+1|2k+1},\nonumber\\
&\Sigma_{2k+2|2k+1}=A\Sigma^\prime_{2k+1|2k}A^T-A\Sigma^\prime_{2k+1|2k}C^T\times\nonumber\\
&\hskip 10mm(C\Sigma^\prime_{2k+1|2k}CT+R+S_{b-r})^{-1}
\Sigma^\prime_{2k+1|2k}A^T+Q.\nonumber
\end{align}
%
\end{lemma}

Although both prediction covariances  $\lim_{k\to \infty} \Sigma_{2k|2k-1}$ and $ \lim_{k\to \infty}\Sigma_{2k+1|2k}$ may have  different limiting values for Strategies II \&\ III,  the value of the former suffices for the rest of the calculation.
\begin{corollary}\label{cor:cor1}
For Strategy I, $\Sigma^{p,\infty}_I\triangleq \lim_{k\to \infty}\Sigma_{k|k-1}$ and  $\Sigma^{\infty}_I\triangleq \lim_{k\to\infty}{\Sigma_{k|k}}$   satisfy
\begin{align}
\label{eq:dare1}
\begin{split}
\Sigma^{p,\infty}_I&=\texttt{dare}(A^T,C^T,Q,R+S_b),\\
\Sigma^{\infty}_I&=\Sigma^{p,\infty}_I-\Sigma^{p,\infty}_IC^T(C\Sigma^{p,\infty}_IC^T+R+S_b)^{-1}C\Sigma^{p,\infty}_I.
\end{split}
\end{align}
\end{corollary}
\begin{corollary}\label{cor:cor2}
For Strategy II,  $\Sigma^{p,\infty}_{II}\triangleq\lim_{k\to \infty}\Sigma_{2k|2k-1}$, $\Sigma^{\infty}_{II_{\text{even}}}\triangleq\lim_{k\to\infty}{\Sigma_{2k|2k}}$ and $\Sigma^{\infty}_{II_{\text{odd}}}\triangleq\lim_{k\to\infty}{\Sigma_{2k+1|2k+1}}$  
 satisfy
\begin{align}
\label{eq:dareII2}
\begin{split}
\Sigma^{p,\infty}_{II}&=\texttt{dare}\left(A^{2^T},C^T, AQA^T+Q,R+S_{2b}\right),\\
\Sigma^{\infty}_{II_{\text{even}}}&=\Sigma^{p,\infty}_{II}-\Sigma^{p,\infty}_{II}C^T(C\Sigma^{p,\infty}_{II}C^T+R+S_b)^{-1}C\Sigma^{p,\infty}_{II},\\
\Sigma^{\infty}_{II_{\text{odd}}}&=\Sigma^{p,\infty}_{II}-\Sigma^{p,\infty}_{II}C^T(C\Sigma^{p,\infty}_{II}C^T+R+S_{2b})^{-1}C\Sigma^{p,\infty}_{II}.
\end{split}
\end{align}
\end{corollary}
\begin{small}
\begin{corollary}\label{cor:cor3}
For Strategy III,  $\Sigma^{p,\infty}_{III}\lim_{k\to \infty}\triangleq\Sigma_{2k|2k-1}$,  
$\Sigma^{\infty}_{III_{\text{even}}}\triangleq\lim_{k\to\infty}{\Sigma_{2k|2k}},$ and $\Sigma^{\infty}_{III_{\text{odd}}}\triangleq\lim_{k\to\infty}{\Sigma_{2k+1|2k+1}},$   
 satisfy
\begin{align}
\label{eq:dareII}
\begin{split}
\Sigma^{p,\infty}_{III}&=\texttt{dare}\left(A^{2^T},\mathcal G_1,AQA^T+Q,\mathcal G_2,\begin{bmatrix}0&AQC^T\end{bmatrix},\texttt{eye}(n)\right),
\\
\Sigma^{\infty}_{III_{\text{even}}}&=\Sigma^{p,\infty}_{III}-\Sigma^{p,\infty}_{III}C^T(C\Sigma^{p,\infty}_{III}C^T+R+S_{b+r})^{-1}C\Sigma^{p,\infty}_{III},\\
\Sigma^{\infty}_{III_{\text{odd}}}&=\Sigma^{p,\infty}_{III}-\Sigma^{p,\infty}_{III}C^T(C\Sigma^{p,\infty}_{III}C^T+R+S_{b-r})^{-1}C\Sigma^{p,\infty}_{III},
\end{split}
\end{align}
\end{corollary}
where
\begin{align*}
\mathcal G_1=\begin{bmatrix}C^T&A^TC^T\end{bmatrix},~~\mathcal G_2=\begin{bmatrix}R+S_{b+r}&0\\0&CQC^T+R+S_{b-r}\end{bmatrix}.
\end{align*}
\end{small}

\section{Control performance analysis}
The limiting performance of three strategies may be computed using standard covariance methods.
\begin{definition}
The  $i,j$-block $(n\times n)$ entry of matrices $\Psi_{X}$, below is denoted by $\Psi_{X}(i,j)$ for $X=I,II,III$.
\end{definition}
\begin{theorem}\label{the:the1}
Subject to Assumption~\ref{ass:nosat}, the performance for Strategy I given by
\begin{align}
J_I=\text{trace}[Q_c\Psi_I(1,1)]+\text{trace}[K^TR_cK\Psi_I(2,2)],\label{eq:perf1}
\end{align}
calculated through these steps:
\begin{enumerate}[label=(\roman*)]
\item\label{step:ARE} $K=\texttt{dare}(A,B,Q_c,R_c),$
\item $\Sigma^{p,\infty}_I=\texttt{dare}\left(A^T,C^T,Q,R+S_b\right),\label{eq:sigma1}$
\item $L=\Sigma^{p,\infty}_I C^T\left(C\Sigma^{p,\infty}_I C^T+R+S_b\right)^{-1},\label{eq:l1}$
\item $\Psi_I=\texttt{dlyap}\left(\mathcal M_1, \mathcal N_1\mathcal P_1\mathcal N^T_1\right),$
\end{enumerate}
where
\begin{align*}
\mathcal M_1&=\begin{bmatrix}A&-BK\\LCA&(I-LC)A-BK\end{bmatrix},~~\mathcal N_1=\begin{bmatrix}I&0&0\\LC&L&L\end{bmatrix},\end{align*}
\begin{align*}
\mathcal P_1=\begin{bmatrix}Q&0&0\\0&R&0\\0&0&S_b\end{bmatrix},~~\Psi_I&=\begin{bmatrix}E(x_k x_k^T)&E(x_k\hat x_{k|k}^T)\\E(\hat x_{k|k}x_k^T)&E(\hat x_{k|k}\hat x_{k|k}^T)\end{bmatrix}.
\end{align*}
\end{theorem}
\begin{theorem} \label{the:the2}
Subject to Assumption~\ref{ass:nosat}, the performance for Strategy II, given by 
\begin{align}
J_{II}&=\frac{1}{2}\text{trace}\left\{Q_c[\Psi_{II}(1,1)+\Psi_{II}(3,3)]\right\}+\nonumber\\     
    &\hskip 10mm  \frac{1}{2}\text{trace}\left\{K^TR_cK[\Psi_{II}(2,2)+\Psi_{II}(4,4)]\right\},\label{eq:perf2}
\end{align}
 calculated through these steps:
\begin{small}
\begin{enumerate}[label=(\roman*)]
\item $K=\texttt{dare}(A,B,Q_c,R_c),$
\item $\Sigma^{p,\infty}_{II}=\texttt{dare}\left(A^{2^T},C^T, AQA^T+Q,R+S_{2b}\right),\label{eq:dare2}$
\item $L_{\text{even}}=\Sigma^{p,\infty}_{II} C^T\left(C\Sigma^{p,\infty}_{II} C^T+R+S_b\right)^{-1}.\label{eq:l1}$
\item $L_{\text{odd}}=\Sigma^{p,\infty}_{II}C^T(C\Sigma^{p,\infty}_{II}C^T+R+S_{2b})^{-1},\label{eq:l2kp1}$
\item $\Psi_{II}=\texttt{dlyap}\left(\mathcal M_2, \mathcal N_2\mathcal P_2\mathcal N_2^T\right)$
\end{enumerate}
\end{small}
where
\begin{small}
\begin{align*}
~~~\mathcal M_2&=F_4F_3F_2F_1,~~\mathcal N_2=\begin{bmatrix}F_4F_3F_2G_1&F_4G_3&F_4F_3G_2&G_4\end{bmatrix},\\
\mathcal P_2&=\begin{bmatrix}Q&0&0&0\\0&Q&0&0\\0&0&R+S_{b}&R+S_{2b}\\0&0&R+S_{2b}&R+S_{2b}\end{bmatrix},\\
\Psi_{II}&=E\left(\begin{bmatrix}x_{2k}\\\hat x_{2k|2k}\\x_{2k+1}\\\hat x_{2k+1|2k+1}\end{bmatrix}\begin{bmatrix}x^T_{2k}&\hat x^T_{2k|2k}&x^T_{2k+1}&\hat x^T_{2k+1|2k+1}\end{bmatrix}\right)\\
G_1&=\begin{bmatrix}0\\I\end{bmatrix},~~G_2=\begin{bmatrix}0\\0\\L_{0}\end{bmatrix},~~G_3=\begin{bmatrix}0\\0\\L_{1}\end{bmatrix},~~G_4=\begin{bmatrix}0\\0\\I\\0\end{bmatrix}\\
F_1&=\begin{bmatrix}0&0&0&A-BK\\0&0&A&-BK\end{bmatrix},~~F_2=\begin{bmatrix}I&0\\0&I\\(I-L_{0}C)&L_{0}C\end{bmatrix},\\
F_3&=\begin{bmatrix}0&I&0\\0&0&I\\(I-L_{1}C)&L_{1}C&0\end{bmatrix},~~F_4=\begin{bmatrix}I&0&0\\0&I&0\\A&-BK&0\\0&-BK&A\end{bmatrix}.
\end{align*}
\end{small}
\end{theorem}
Note, the two-step update is described by the recursion
\begin{align*}
\begin{bmatrix}x_{2k}\\\hat x_{2k|2k}\\x_{2k+1}\\\hat x_{2k+1|2k+1}\end{bmatrix}
&=\mathcal M_2\begin{bmatrix}x_{2k-2}\\\hat x_{2k-2|2k-2}\\x_{2k-1}\\\hat x_{2k-1|2k-1}\end{bmatrix}
+\mathcal N_2\begin{bmatrix}w_{2k-1}\\w_{2k}\\\begin{bmatrix}v_{2k}+q_{2k}\\v_{2k}+q_{2k+1}\end{bmatrix}
\end{bmatrix}.
\end{align*}

\begin{theorem}\label{the:the3}
Subject to Assumption~\ref{ass:nosat}, the performance for Strategy III, given by
\begin{align}
J_{III}&=\frac{1}{2}\text{trace}\left\{Q_c[\Psi_{III}(1,1)+\Psi_{III}(3,3)]\right\}+\nonumber\\     
    &\hskip 10mm  \frac{1}{2}\text{trace}\left\{K^TR_cK[\Psi_{III}(2,2)+\Psi_{III}(4,4)]\right\},\label{eq:perf3}
\end{align}
 calculated through these steps:
\begin{small}
\begin{enumerate}[label=(\roman*)]
\item $K=\texttt{dare}(A,B,Q_c,R_c),$
\item $\Sigma^{p,\infty}=\texttt{dare}\left(A^{2^T},\mathcal G_1,AQA^T+Q,\mathcal G_2,\begin{bmatrix}0&AQC^T\end{bmatrix},\texttt{eye}(n)\right),$
\item $L_{\text{even}}=\Sigma^{p,\infty} C^T\left(C\Sigma^{p,\infty} C^T+R+S_b\right)^{-1},\label{eq:l30}$
\item $L_{\text{odd1}}=\Sigma^{p,\infty}C^T(C\Sigma^{p,\infty}C^T+R+S_{b+r})^{-1},\label{eq:l31}$
\item $L_{\text{odd2}}=\Sigma^{p,\infty}C^T(C\Sigma^{p,\infty}C^T+R+S_{b-r})^{-1},\label{eq:l32}$
\item $\Psi_{III}=\texttt{dlyap}\left(\mathcal M_3, \mathcal N_3\mathcal P_3\mathcal N_3^T\right),$
\end{enumerate}
\end{small}
where
\begin{small}
\begin{align*}
\mathcal G_1&=[C^T A^TC^T],\\
\mathcal G_2&=\begin{bmatrix}R+S_{b+r}&0\\0&CQC^T+R+S_{b-r}\end{bmatrix}\\
~~~\mathcal M_2&=F_4F_3F_2F_1,\\
~~\mathcal N_2&=\begin{bmatrix}F_4F_3F_2G_1&F_4G_3&F_4F_3G_2&G_4\end{bmatrix},\\
\mathcal P_3&=\begin{bmatrix}Q&0&0&0&0\\0&Q&0&0&0\\0&0&R+S_b&R+S_{b+r}&0\\0&0&R+S_{b+r}&R+S_{b+r}&0\\0&0&0&0&R+S_{b-r}\end{bmatrix}\\
\Psi_{III}&=E\left(\begin{bmatrix}x_{2k+1}\\\hat x_{2k+1|2k+1}\\x_{2k}\\\hat x^1_{2k|2k}\end{bmatrix}
\begin{bmatrix}x^T_{2k+1}&\hat x^T_{2k+1|2k+1}&x^T_{2k}&\hat x^{1^T}_{2k|2k}\end{bmatrix}\right)\\\\
G_1&=\begin{bmatrix}I\\0\end{bmatrix},~~G_2=\begin{bmatrix}0&0\\L_{\text{even}}&0\\0&L_{\text{odd1}}\end{bmatrix},~~G_3=\begin{bmatrix}I\\0\\0\\0\end{bmatrix},\\
~~G_4&=\begin{bmatrix}0\\L_{\text{odd2}}\\0\\0\end{bmatrix},~~F_1=\begin{bmatrix}A&-BK&0&0\\0&A-BK&0&0\end{bmatrix},\\
F_3&=\begin{bmatrix}A&-BK&0\\0&0&I\\I&0&0\\0&I&0\end{bmatrix}, ~~F_2=\begin{bmatrix}I&0\\L_{\text{even}}C&(I-L_{\text{even}}C)\\L_{\text{odd1}}C&(I-L_{\text{odd1}}C)\end{bmatrix},\\
F_4&=\begin{bmatrix}I&0&0&0\\L_{\text{odd2}}C&(I-L_{\text{odd2}}C)A&0&-(I-L_{\text{odd2}}C)BK\\
0&0&I&0\\0&0&0&I\end{bmatrix},
\end{align*}
\end{small}
\end{theorem}
The two-step update is described by the recursion
\begin{align*}
\begin{bmatrix}x_{2k+1}\\\hat x_{2k+1|2k+1}\\x_{2k}\\\hat x^1_{2k|2k}\end{bmatrix}
&=\mathcal M_3\begin{bmatrix}x_{2k-1}\\\hat x_{2k-1|2k-1}\\x_{2k-1}\\\hat x^1_{2k-2|2k-2}\end{bmatrix}
 \nonumber\\&+\mathcal N_3\begin{bmatrix}w_{2k-1}\\w_{2k}\\\begin{bmatrix}v_{2k}+q_{b,2k}\\v_{2k}+q_{b+r,2k}\end{bmatrix}\\v_{2k+1}+q_{b-r,2k+1}\end{bmatrix}.
\end{align*}

\section{Escape time analysis}\label{sec:esc}
The performance analysis from earlier sections is based on direct
second moment calculations subject to the validity of
Assumption~\ref{ass:nosat}, i.e. that the controlled system dithered output 
\begin{align*}
z_t=y_t+d_t,
\end{align*}
satisfies $|z_t|\leq\zeta$. For
Gaussian $y_t$, or indeed for any $y_t$ with density of unbounded support, the signal $y_t+d_t$ is guaranteed to exceed this bound infinitely often. Our aim in this section is to quantify the average residence time of the dithered controlled output signal inside the saturation bound. If this residence time is long, then the earlier linear analysis will remain valid on average for a long time and can still be used to characterize performance, since the stabilizing control yields a quasi-stationary closed loop subject to no saturation. This will be validated by computational experiments in Section~\ref{sec:numEx}.

We make the following definition.
\begin{definition}
The \textit{escape time,} $\tau_\text{esc},$ is the first time that $z_t\not\in[-\zeta,\zeta]$.
\end{definition}
Our aim is now to calculate the mean escape time as a function of
$\zeta$. This will demonstrate that the choice of $\zeta$ to yield a
particular mean escape time depends on the choice of state feedback
control gain $K$. The state estimation covariance analysis of
Section~\ref{sec:KF} did not explicitly depend on $K$. But now, via
its effect on $\zeta,$ the control problem affects this covariance.

 If we have ergodicity of the stochastic process $\{z_t\}$ then the
 long-term sample average frequency of $z_t$ falling outside
 $[-\zeta,\zeta]$ is equal to the ensemble average computable from the
 density of $z_t.$ If the Gaussian process $\{y_t\}$ is ergodic, then since $\{d_t\}$ is white and stationary, the signal $\{y_t+d_t\}$ is  ergodic. We have the following theorem from Caines \cite{caines2018linear} who cites earlier sources going back to Maruyama and Grenander.
\begin{theorem}\cite{caines2018linear}\label{the:Rozanov}
A necessary and sufficient condition for a discrete-time stationary Gaussian process to be ergodic is that the spectral distribution of the process is continuous.
\end{theorem}

If $y_t$ is the output of a stable linear system  driven by white, independent, zero-mean Gaussian noises $n_t$ and $r_t$ with covariances $Q$ and $R$ respectively, that is,
\begin{align*}
\xi_{t+1}&=F\xi_t+Gn_t,\\
y_t&=Hp_t+Jr_t,
\end{align*}
then its power spectral density is given by
\begin{align*}
\Phi_{yy}(\omega)&=JRJ^T+H(e^{j\omega}I-F)^{-1}GQG^T(e^{-j\omega} I -F^T)^{-1}.
\end{align*}
If the eigenvalues of $F$ are within $|z|<1$ and $JRJ^T>0,$ then $y_t$ is ergodic
by Theorem~\ref{the:Rozanov}, as is the signal $z_t$.

For our LQG problem, $y_t$ is generated with
\begin{align*}
F&=\begin{bmatrix}A&-BK\\LCA&A-BK-LCA\end{bmatrix},~~~ G=\begin{bmatrix}I&0\\LC&L\end{bmatrix},\\
H&=\begin{bmatrix}C&0\end{bmatrix},~~~ J=I,
\end{align*}
which $F$ has all eigenvalues inside the unit circle by construction subject to the conditions in the following theorem.

\begin{theorem}\cite{soderstrom2012discrete}\label{the:ergodic}
Subject to Assumption~\ref{ass:nosat}, provided $R_c>0$, $R>0$, $[A,Q_c]$ detectable, $[A,Q]$ stabilizable, the dithered controlled output signal, $z_t=y_t+d_t$, is asymptotically stationary and ergodic. So
\begin{align}
\lim_{T\to\infty}{\frac{1}{T}\sum_{t=1}^T{\mathbb 1_{|z_t|>\zeta}}}&=\text{Pr}\left(|z_t|>\zeta\right),\label{eq:argo}
\end{align}
where $\mathbb 1_A$ is the indicator function of event $A$.
\end{theorem}

Once we have ergodicity of the closed-loop signal $z_t$, then we have the following result.
\begin{theorem}\label{the:esc}
For ergodic $z_t,$ if Pr$(|z_k|>\zeta)=\beta,$ then the expected escape time is $\text{E}[\tau_{\text{esc}}]=\frac{1}{\beta}.$
\end{theorem} 

 These are the steps and important parameters of the analysis.
 \begin{enumerate}
 \item Choose a desired average escape time $\text{E}[\tau_\text{esc}].$ The escape probability is $\beta=\frac{1}{\text{E}[\tau_\text{esc}]}.$
 \item Initiate the line search for $\zeta$ 
depending on coding strategy and using one of
 \begin{itemize}
 \item $\Psi _I(1,1)$ from \eqref{eq:perf1}, or
 \item $\Psi_{II}(1,1)$ and $\Psi_{II}(3,3)$ from \eqref{eq:perf2}, or
 \item $\Psi_{III}(1,1)$ and $\Psi_{III}(3,3)$ from \eqref{eq:perf3},
 \end{itemize}
compute $Z$, the covariance of $z_t$. Then solve
 \begin{align*}
 \frac{\beta}{2}&=\texttt{mvncdf}\left(-\zeta_\text{new}.\mathbb 1_m,0_m,Z\right),
 \end{align*}
for $\zeta_\text{new}$ where \texttt{mvncdf} is the multivariate normal cumulative distribution function.
 \end{enumerate}
\section{Numerical examples\label{sec:numEx}}
We compare  coding strategies in  the   following examples through these steps. 
\begin{enumerate}
\item With given $\{A,B,C,Q, R,Q_c, R_c\}$, compute the linear feedback gain $K,$ via Theorem~\ref{the:the1} Step~\ref{step:ARE}.
\item Fix the mean escape time, $\tau_\text{esc}.$
\item For each coding strategy, compute the corresponding quantizer bound, $\zeta,$ using the iteration described below Theorem~\ref{the:esc}.
\item Compute the performance of each strategy using Theorems~\ref{the:the1}-\ref{the:the3}, as appropriate.  
\end{enumerate}

\subsection{Escape time and  quantizer bound}\label{B}

We compute the escape time through two methods, the analytical method based on Theorem  \ref{the:esc} and  simulation. In addition, we compare the performance of the coding strategies. Let us define the parameters as follows. 
\begin{itemize}
\item $R_c,$ control weight in LQ output feedback control.
\item $A-BK,$ LQ closed-loop matrix.
\item $\zeta,$ quantization bound.
\item $\tau_a,$ mean escape time computed via Theorem \ref{the:esc}.
\item $\tau_\text{emp},$ empirical mean escape time from simulation.
\item $J_I, J_{II},$ corresponding performances for Strategies I and II.
\end{itemize}
In  the  simulation for computing, $\tau_\text{emp},$ we  average over 20,000 iterations the first time that the dithered output signal  breaches the quantizer bound for  a 5,000-sample simulation with the following parameters for a scalar system. Then we compare the performance of two different strategies with a fixed time $\tau=1000$ and   parameters as follows.
\begin{verbatim}
A = 0.9999; B = 1; C = 1; Q = 1; 
R = 1; Qc = 1; 
\end{verbatim}  
for 3-bit quantizer
\begin{small}
\begin{center}
\begin{tabular}{|l|c|c|c|c|c|c|}
\hline
$R_c$&A-BK&$\zeta$&$\tau_a$&$\tau_{\text{emp}}$&$J_I$&$J_{II}$\\
\hline
 1e5&0.9968&43.14&1000&2320&325&309\\
 1e4&0.9900&24.67&1000&2194&104&101\\
 1e3&0.9689&14.50&1000&1813&34.136&34.135\\
 100&0.9049&9.15&1000&1317&11.81&12.43\\
 10&0.7298&6.62&1000&1040&4.78&5.56\\
 1&0.3819&5.68&1000&990&2.64&3.45\\  
 0.1&0.0839&5.49&998&977&2.10&2.92\\   
\hline
\end{tabular}
\end{center}
\end{small}  
 for 2-bit quantizer
\begin{small}
\begin{center}
\begin{tabular}{|l|c|c|c|c|c|c|}
\hline
$R_c$&A-BK&$\zeta$&$\tau_a$&$\tau_{\text{emp}}$&$J_I$&$J_{II}$\\
\hline
 1e5&0.9968&48.14&1000&2354&474&315\\
 1e4&0.9900&27.74&1000&2159&137&103\\
 1e3&0.9689&16.44&1000&1780&42.22&35.04\\
 100&0.9049&10.43&1000&1413&14.34&12.97\\
 10&0.7298&7.54&1000&1213&5.94&5.91\\
 1&0.3819&6.40&1000&1164&3.43&3.73\\ 
 0.1&0.0839&6.12&998&1146&2.81&3.18\\  
\hline
\end{tabular}
\end{center}
\end{small}  
\begin{figure}[h]
\begin{center}
\includegraphics[trim={0 7cm 0 5cm},clip,scale=.45]{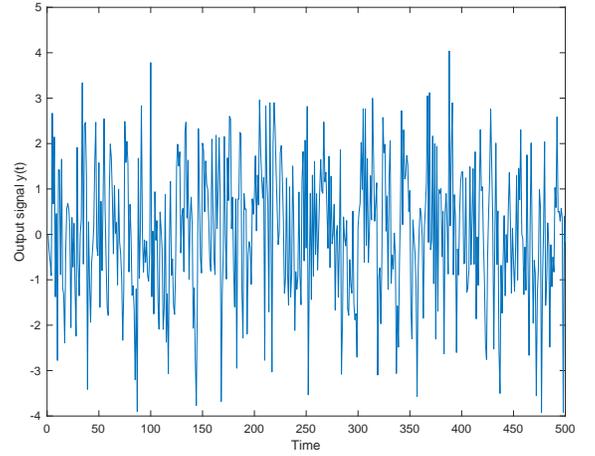}
\caption{Controlled output signal y(t) with 3-bit coding and $R_c=0.01$, corresponding to roughly minimum-variance control and hence to low amplitude, near-white $y_t$. Coding provides little benefit.\label{unpred}}
\end{center}
\end{figure}
\begin{figure}[h]
\begin{center}
\includegraphics[trim={0 7cm 0 5cm},clip,scale=.45]{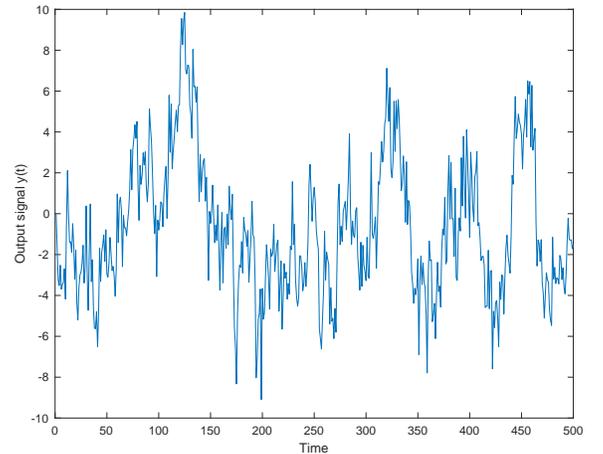}
\caption{Controlled output signal y(t) with 3-bit coding and $R_c=100$, corresponding to higher amplitude, correlated $y_t$. Coding provides tangible control benefit.\label{pred}}
\end{center}
\end{figure}


%
%

%

As we may conclude from the above example, the coding strategy is picked based on the nature of the controlled output signal. If the output signal has a random or unpredictable nature, as in Figure \ref{unpred},  the coding has less benefits and we stick with Strategy I. 
 In contrast, the coding strategy has advantages if the output controlled signal is more regulated or predictable such as Figure \ref{pred}. In this case, as we have higher resolution or accuracy by including  a delay in updating the measurement,  Strategy II outperforms Strategy I in which  the measurement is updated each time but with less accuracy.
 When the control objective is minimum variance, the output signal resembles to a white noise signal and the quantization bound has  smaller size, so the coding has no benefits. Once we move away from a minimum variance control objective, the output signal $y_k$  has both larger amplitude and an output signal $y_k$ which is more correlated. So coding can bring tangible performance benefits.
 
The computational examples exhibit the following.
\begin{itemize}
\item The control performance advantage of coding is tied to the redundancy in the regulated system output. 
\item  Delay-one minimum variance control benefits little from coding, because the regulated output is close to white.
\item If the number of bits, $b,$ is large then coding has limited benefit.
\item For a given escape time, the quantizer bound, $ \zeta,$ is smaller for  better regulated signal $y_t.$
\end{itemize}
\section{Conclusion}
We have explored three very specific periodic coding strategies of the plant output signal and their effect on LQ performance subject to an expected escape time. The interaction between the control law and the estimation problem occurs through the selection of the upper bound, $\zeta,$ of the dithered quantizers.  The general conclusion is that the more correlated is the controlled output, the more benefit is achieved by coding. So that minimum variance problems should exhibit less gain from coding than should those with heavier control penalty. The computational examples show that these coding schemes promise most value when the number of bits is small. These are generalizable conclusions to other more sophisticated codes and reflect observations in signal processing, but without the connection to a control objective.

The novelties of the approach lie in the treatment of  dithered quantizers and the introduction of the system escape time as a tool for analysis. The decomposition of the quantizer into two parts -- infinite quantizer plus saturation -- together with the escape time permits the consideration of linear controlled covariances and the distinct escape time analysis. This study of escape time is distinguished from many other studies which seek to manage asymptotic or infinite-horizon average or moment properties.

\section*{Appendix}
\emph{Proof for lemma \ref{lem:lem2}:}\\
Let us start with initial state estimate $\hat x_{0|-1}$ and covariance $\Sigma_{0|-1},$ the  Kalman filter is calculated by:
\vskip 1mm
\noindent\textit{At even time, $t=2k:$ (Low resolution measurement)}
\begin{align}
&L_{2k}=\Sigma_{2k|2k-1}C^T(C\Sigma_{2k|2k-1}C^T+R+S_b)^{-1},\nonumber\\
&\hat x_{2k|2k}= \hat x_{2k|2k-1}+L_{2k}(p_{2k}-C\hat x_{2k|2k-1}),\nonumber
\end{align}
\vskip 1mm
\noindent\textit{At odd time, $t=2k+1:$ (High resolution measurement)}
We receive the less significant part of the quantized $y_{2k}$ and construct the 2$b$-bits measurement $z_{2b,2k+1}=z_{b,2k}\oplus_{2b,2k}$ through concatenation,
\begin{align}
&p_{2k+1}\gets p_{2k}+2^{-b}m_{2k+1}-d^{2b}_{2k}\nonumber\\
&L_{2k+1}=\Sigma_{2k|2k-1}C^T(C\Sigma_{2k|2k-1}C^T+R+S_{2b})^{-1},\nonumber\\
&\hat x_{2k+1|2k+1}=A\hat x_{2k|2k+1}+Bu_{2k},\nonumber\\
&\hat x_{2k+1|2k+1}= A\left(\hat x_{2k|2k-1}+L_{2k+1}(p_{2k+1}-C\hat x_{2k|2k-1})\right)\nonumber\\
&\hskip 58mm-BK\hat x_{2k|2k},\nonumber\\
&\hat x_{2k+2|2k+1}=(A-BK)\hat x_{2k+1|2k+1},\nonumber\\
&\Sigma_{2k+2|2k+1}=A^2\Sigma_{2k|2k-1}{A^2}^T-A^2L_{2k+1}C\Sigma_{2k|2k-1}{A^2}^T\nonumber\\
&\hskip 60mm+AQA^T+Q.\nonumber
\end{align}
\emph{Proof for Lemma \ref{lem:lem3}:}\\
Let us start with initial state estimate $\hat x_{0|-1}$ and covariance $\Sigma_{0|-1}, \text{at } t=0,$ the  Kalman filter is calculated by:
\vskip 1mm\noindent\textit{At even time, $t=2k,$}
\begin{align}
&L_{2k}=\Sigma_{2k|2k-1}C^T(C\Sigma_{2k|2k-1}C^T+R+S_b)^{-1},\nonumber\\
&\hat x_{2k|2k}= \hat x_{2k|2k-1}+L_{2k}(p_{2k}-C\hat x_{2k|2k-1}),\nonumber
\end{align}
\vskip 1mm\noindent
\textit{At odd times, $t=2k+1,$}
\begin{align}
&p^\prime_{2k}=p_{2k}+2^{-b}\text{MSB}_r(m_{2k+1})-d^{b+r}_{2k},\nonumber\\
&L^\prime_{2k}=\Sigma_{2k|2k-1}C^T(C\Sigma_{2k|2k-1}C^T+R+S_{b+r})^{-1},\nonumber\\
&\hat x^\prime_{2k+1|2k}=A\hat x_{2k|2k-1}+AL_{2k}^\prime(p^\prime_{2k}-C\hat x_{2k|2k-1})\nonumber\\
&\hskip 50mm-BK\hat x_{2k|2k},\nonumber\\
&\Sigma^\prime_{2k+1|2k}=A\Sigma_{2k|2k-1}A^T-A\Sigma_{2k|2k-1}C^T\times\nonumber\\
&\hskip 10mm(C\Sigma_{2k|2k-1}C^T+R+S_{b+r})^{-1}C\Sigma_{2k|2k-1}A^T+Q,\nonumber\\
&p_{2k+1}=\text{LSB}_{b-r}(m_{2k+1})-d^{b-r}_{2k+1},\nonumber\\
&L_{2k+1}=\Sigma^\prime_{2k+1|2k}C^T(C\Sigma^\prime_{2k+1|2k}C^T+R+S_{b-r})^{-1},\nonumber\\
&\hat x_{2k+1|2k+1}=\hat x^\prime_{2k+1|2k}+L_{2k+1}(p_{2k+1}-C\hat x^\prime_{2k+1|2k}),\nonumber\\
&\hat x_{2k+2|2k+1}=A\hat x_{2k+1|2k+1}+Bu_{2k+1},\nonumber\\
&\hat x_{2k+2|2k+1}=(A-BK)\hat x_{2k+1|2k+1},\nonumber\\
&\Sigma_{2k+2|2k+1}=A\Sigma^\prime_{2k+1|2k}A^T-A\Sigma^\prime_{2k+1|2k}C^T\times\nonumber\\
&\hskip 10mm(C\Sigma^\prime_{2k+1|2k}CT+R+S_{b-r})^{-1}
\Sigma^\prime_{2k+1|2k}A^T+Q.\nonumber
\end{align}

\emph{Proof for Corollary \ref{cor:cor1}:}\\
Can be found at 
\cite{anderson1979optimal} but the difference is quantization noise $S_b$ is added to the measurement noise, and it is replaced by $R+S_b$ in all calculation.\\
\emph{Proof for Corollary \ref{cor:cor2}:}\\
Let us start with $\Sigma\triangleq\Sigma_{2k|2k-1}$
\begin{enumerate}[label=(\roman*)]
\item Low resolution measurement 
\begin{align}
\Sigma_{2k|2k}=\Sigma -\Sigma C^T(C\Sigma C^T+R+S_b)^{-1}C\Sigma \label{equ:covf}
\end{align}
\item High resolution measurement
\begin{align}
\Sigma_{2k|2k+1}=\Sigma-\Sigma C^T(C\Sigma C^T+R+S_{2b})^{-1}C\Sigma \label{equ:covs}
\end{align}
\item Time $2k+1$ filtered measurement
\begin{align*}
\Sigma_{2k+1|2k+1}&=A\Sigma_{2k|2k+1}A^T+Q.
\end{align*}
\item Time $2k+2$ prediction and let $k\to\infty$
\begin{align*}
\Sigma_{2k+2|2k+1}&=A\Sigma_{2k+1|2k+1}A^T+Q,\\
&=A^2\Sigma_{2k|2k+1}A^{2^T}+AQA^T+Q,
\end{align*}
\begin{align*}
\Sigma^{p,\infty}_{II}&=A^2\Sigma^{p,\infty}_{II} A^{2^T}-A^2\Sigma^{p,\infty}_{II} C^T\times \nonumber\\
&\hskip 10mm(C\Sigma^{p,\infty}_{II} C^T+R+S_{2b})^{-1}C\Sigma^{p,\infty}_{II} A^{2^T}+\\
& \hskip 35mm AQA^T+Q.
\end{align*}
where
\begin{small}
\begin{align}
\lim_{k\to \infty}\Sigma_{2k|2k-1}=\Sigma^{p,\infty}_{II}=\texttt{dare}\left(A^{2^T},C^T, AQA^T+Q,R+S_{2b}\right).\nonumber
\end{align}
\end{small}
\item Substitute $\Sigma^{p,\infty}_{II}$ into \eqref{equ:covf} and \eqref{equ:covs} 
\begin{small}
\begin{align*}
\Sigma^{p,\infty}_{II}&=\texttt{dare}\left(A^{2^T},C^T, AQA^T+Q,R+S_{2b}\right),\\
\Sigma^{\infty}_{II_{\text{even}}}&=\Sigma^{p,\infty}_{II}-\Sigma^{p,\infty}_{II}C^T(C\Sigma^{p,\infty}_{II}C^T+R+S_b)^{-1}C\Sigma^{p,\infty}_{II},\\
\Sigma^{\infty}_{II_{\text{odd}}}&=\Sigma^{p,\infty}_{II}-\Sigma^{p,\infty}_{II}C^T(C\Sigma^{p,\infty}_{II}C^T+R+S_{2b})^{-1}C\Sigma^{p,\infty}_{II}.
\end{align*}
\end{small}
\end{enumerate}
\emph{Proof for Corollary \ref{cor:cor3}:}\\

The period-two update consists of two pieces starting from the same initial data, $(\hat x_{2k|2k-1},\Sigma=\Sigma_{2k|2k-1})$.
\begin{description}
\item[Even times] -- No need to keep track of this in the computation of the covariance $\Sigma_{2k+1|2k+1}$ since this is calculated based only on 
\begin{align*}
z^\prime_{2k}&=Cx_{2k}+v_{2k}+q_{b+r,2k},\\
z_{2k+1}&=Cx_{2k+1}+v_{2k+1}+q_{b-r,2k+1}.
\end{align*}
It is, however, important for the Lyapunov computation.
\item[Odd times] -- We skip over the even step and use both $z^\prime_{2k}$ and $z_{2k+1}$ to update $\hat x_{2k|2k-1}.$ Start as usual.
\begin{align*}
x_{2k+1}&=Ax_{2k}+Bu_{2k}+w_{2k},\\
p_{2k+1}&=Ax_{2k}+w_{2k},\\
z^\prime_{2k}&=Cx_{2k}+v_{2k}+q_{b+r,2k},\\
z_{2k+1}&=CAx_{2k}+CBu_{2k}+Cw_{2k}+v_{2k+1}+q_{b-r,2k+1},\\
\zeta_{2k+1}&=CAx_{2k}+Cw_{2k}+v_{2k+1}+q_{b-r,2k+1}
\end{align*}
where denote $\zeta_{2k+1}=z_{2k+1}-CBu_{2k}$ and calculate joint conditional density,
\begin{align*}
\text{pdf}\left(\left.\begin{bmatrix}p_{2k+1}\\z^\prime_{2k}\\\zeta_{2k+1}\end{bmatrix}\right|\mathcal Z^{2k-1}\right)&=
\mathcal N\left(\begin{bmatrix}A\hat x_{2k|2k-1}\\C\hat x_{2k|2k-1}\\CA\hat x_{2k|2k-1}\end{bmatrix},
\mathcal M\right),
\end{align*}
\begin{tiny}
\begin{align*}
\mathcal M&=\begin{bmatrix}A\Sigma A^T+Q&A\Sigma C^T&A\Sigma A^TC^T+QC^T\\
C\Sigma A^T&C\Sigma C^T+R+S_{b+r}&C\Sigma A^TC^T\\
CA\Sigma A^T+CQ&CA\Sigma C^T&CA\Sigma A^TC^T+CQC^T+R+S_{b-r}\end{bmatrix},
\end{align*}
\end{tiny}
hence,
\begin{small}
\begin{align*}
&\text{cov}(x_{2k+1}\vert\mathcal Z^{2k+1})=A\Sigma A^T+Q-\begin{bmatrix}A\Sigma C^T&A\Sigma A^TC^T+QC^T\end{bmatrix}\times
\\
&\hskip 5mm\begin{bmatrix}C\Sigma C^T+R+S_{b+r}&C\Sigma A^TC^T\\CA\Sigma C^T&CA\Sigma A^TC^T+CQC^T+R+S_{b-r}\end{bmatrix}
^{-1}\times\nonumber\\
&\hskip 15mm\begin{bmatrix}C\Sigma A^T\\CA\Sigma A^T+CQ\end{bmatrix},
\end{align*}
\end{small}
by taking limits as 
\begin{align*}
\lim_{k\to \infty}\Sigma_{2k|2k-1}=\Sigma_{2k+2|2k+1}=\Sigma,
\end{align*}
\begin{align*}
\text{cov}(x_{2k+2}\vert\mathcal Z^{2k+1})&=A\times\text{cov}(x_{2k+1}\vert\mathcal Z^{2k+1})\times A^T+Q,
\end{align*}
\begin{small}
\begin{align}
\Sigma&=A^2\Sigma A^{2^T}+AQA^T+Q-\begin{bmatrix}A^2\Sigma C^T&A^2\Sigma A^TC^T+AQC^T\end{bmatrix}\times\nonumber\\
&\hskip 2mm\begin{bmatrix}C\Sigma C^T+R+S_{b+r}&C\Sigma A^TC^T\\CA\Sigma C^T&CA\Sigma A^TC^T+CQC^T+R+S_{b-r}\end{bmatrix}
^{-1}\times \nonumber\\
&\hskip 25mm\begin{bmatrix}C\Sigma A^{2^T}\\CA\Sigma A^{2^T}+CQA^T\end{bmatrix},\nonumber\\
 &=A^2\Sigma A^{2^T}+AQA^T+Q-\left(A^2\Sigma\begin{bmatrix}C^T&A^TC^T\end{bmatrix}
+\begin{bmatrix}0&AQC^T\end{bmatrix}\right)\times\nonumber \\
&\hskip -2mm\left(\begin{bmatrix}C\\CA\end{bmatrix}\Sigma\begin{bmatrix}C^T&A^TC^T\end{bmatrix}+
\begin{bmatrix}R+S_{b+r}&0\\0&CQC^T+R+S_{b-r}\end{bmatrix}\right)^{-1}\times\nonumber\\
&\hskip 15mm\left(\begin{bmatrix}C\\CA\end{bmatrix}\Sigma A^{2^T}+\begin{bmatrix}0\\CQA^T\end{bmatrix}\right).\nonumber
\end{align}
\end{small}
and we use DARE to calculate,
\begin{tiny}
\begin{align}\label{eq:dare3}
\Sigma^{p,\infty}_{III}&=\texttt{dare}\left(A^{2^T},\begin{bmatrix}C^T&A^TC^T\end{bmatrix},AQA^T+Q,\mathcal G_2,\begin{bmatrix}0&AQC^T\end{bmatrix},\texttt{eye}(n)\right),\nonumber
\end{align}
\end{tiny}
\end{description}
and similar to proof of Corollary \ref{cor:cor2}
\begin{align*}
\Sigma^{\infty}_{III_{\text{even}}}&=\Sigma^{p,\infty}_{III}-\Sigma^{p,\infty}_{III}C^T(C\Sigma^{p,\infty}_{III}C^T+R+S_{b+r})^{-1}C\Sigma^{p,\infty}_{III},\\
\Sigma^{\infty}_{III_{\text{odd}}}&=\Sigma^{p,\infty}_{III}-\Sigma^{p,\infty}_{III}C^T(C\Sigma^{p,\infty}_{III}C^T+R+S_{b-r})^{-1}C\Sigma^{p,\infty}_{III},
\end{align*}
where,
\begin{align*}
\mathcal G_2=\begin{bmatrix}R+S_{b+r}&0\\0&CQC^T+R+S_{b-r}\end{bmatrix}.
\end{align*}
\emph{Proof for Theorem \ref{the:the1}:}\\
Truncate every sample to $b$ bits, transmit 
\begin{align*}
z_t=Cx_t+v_t+q_{b,t}.
\end{align*}
Kalman filter is stationary and satisfies
\begin{align*}
\Sigma&=\texttt{dare}\left(A^T,C^T,Q,R+S_b\right),\\
L&=\Sigma-\Sigma C^T\left(C\Sigma C^T+R+S_b\right).
\end{align*}
Closed-loop equations
\begin{align}
x_{t+1}&=Ax_t-BK\hat x_{t|t}+w_t,\nonumber\\
\hat x_{t+1|t+1}&=\hat x_{t+1|t}+L(Cx_{t+1}+v_{t+1}+q_{b,t+1}\nonumber\\
&\hskip 25mm -C\hat x_{t+1|t}),\nonumber\\
&=(A\hat x_{t|t}-BK\hat x_{t|t})\nonumber\\
&\hskip 2mm+L\left[C(Ax_t-BK\hat x_{t|t}+w_t)\right]\nonumber\\
&\hskip 5mm+L\left[v_{t+1}+q_{b,t+1}-C(A-BK)\hat x_{t|t}\right],\nonumber\\
&=[(I-LC)(A-BK)-LCBK]\hat x_{t|t}\nonumber\\
&\hskip 10mm+LCAx_t+LCw_t\nonumber\\
&\hskip 15mm+Lv_{t+1}+Lq_{b,t+1},\nonumber\\
&=LCAx_t+[(I-LC)A-BK]\hat x_{t|t}\nonumber\\
&\hskip 10mm+LCw_t+Lv_{t+1}+Lq_{b,t+1},\nonumber\\
\hskip 20mm \begin{bmatrix}x_{t+1}\\\hat x_{t+1|t+1}\end{bmatrix}
&=\begin{bmatrix}A&-BK\\LCA&(I-LC)A-BK\end{bmatrix}
\begin{bmatrix}x_t\\\hat x_{t|t}\end{bmatrix}\nonumber\\
&\hskip 10mm+\begin{bmatrix}I&0&0\\LC&L&L\end{bmatrix}\begin{bmatrix}w_t\\v_{t+1}\\q_{b,t+1}\end{bmatrix}.\nonumber
\end{align}
%
Let us denote
\begin{small}
\begin{align*}
\mathcal A&=\begin{bmatrix}A&-BK\\LCA&(I-LC)A-BK\end{bmatrix},\hskip 2mm\mathcal B=\begin{bmatrix}I&0&0\\LC&L&L\end{bmatrix},
\end{align*}
\begin{align*}
\hskip 10mm\mathcal Q=\begin{bmatrix}Q&0&0\\0&R&0\\0&0&S_b\end{bmatrix},~~\Psi_I=\texttt{dlyap}\left(\mathcal M_1, \mathcal N_1\mathcal P_1\mathcal N^T_1\right),
\end{align*}
\end{small}
\hskip 10mm hence the performance calculation, 
\begin{align*}
J_I&=\text{trace}[Q_c\Psi_I(1,1)]+\text{trace}[K^TR_cK\Psi_I(2,2)],
\end{align*}
\hskip 30 mm where,
\begin{align*}
\Psi_I=\begin{bmatrix}E(x_kx_k^T)&E(x_k\hat x_{k|k}^T)\\E(\hat x_{k|k}x_k^T)&E(\hat x_{k|k}\hat x_{k|k}^T)\end{bmatrix}.
\end{align*}

\emph{Proof for Theorem \ref{the:the2}:}\\
Truncate $y_t$ to $b$ bits at even times $t$ and then to $2b$ bits at odd times $t$. The quantization variances $S_b$ and $S_{2b}$  respectively.

Start with $x_{2k-1}$, $\hat x_{2k-1|2k-1}$ and $\Sigma_{2k|2k-1}$. Compute
\begin{align*}
L_{2k}&=\Sigma_{2k|2k-1}C^T(C\Sigma_{2k|2k-1}C^T+R+S_b)^{-1},\\
L_{2k+1}&=\Sigma_{2k|2k-1}C^T(C\Sigma_{2k|2k-1}C^T+R+S_{2b})^{-1}.
\end{align*}
State and predictor update
\begin{align*}
x_{2k}&=Ax_{2k-1}-BK\hat x_{2k-1|2k-1}+w_{2k-1}.\\
\hat x_{2k|2k-1}&=(A-BK)\hat x_{2k-1|2k-1},
\end{align*}
so rearrange these equations,
\begin{align*}
\begin{bmatrix}\hat x_{2k|2k-1}\\x_{2k}\end{bmatrix}
&=\begin{bmatrix}0&0&0&A-BK\\0&0&A&-BK\end{bmatrix}
\begin{bmatrix}x_{2k-2}\\\hat x_{2k-2|2k-2}\\x_{2k-1}\\\hat x_{2k-1|2k-1}\end{bmatrix}\\
&+\begin{bmatrix}0\\I\end{bmatrix}w_{2k-1}
\end{align*}

Filter update with low resolution measurement, $z_{b,2k}$.
\begin{align*}
\hat x_{2k|2k}&=\hat x_{2k|2k-1}+L_{2k}(z_{b,2k}-C\hat x_{2k|2k-1}),
\end{align*}
rearrange the equation in matrix form,
\begin{small}
\begin{align*}
\begin{bmatrix}\hat x_{2k|2k-1}\\x_{2k}\\\hat x_{2k|2k}\end{bmatrix}
&=\begin{bmatrix}I&0\\0&I\\(I-L_{2k}C)&L_{2k}C\end{bmatrix}
\begin{bmatrix}\hat x_{2k|2k-1}\\x_{2k}\end{bmatrix}\\
&+\begin{bmatrix}0\\0\\L_{2k}\end{bmatrix}v_{2k}
+\begin{bmatrix}0\\0\\L_{2k}\end{bmatrix}q_{2k}.
\end{align*}
\end{small}
\hskip 10mm Filter update with high resolution measurement, 

\hskip 10mm $z_{2k+1}=z_{b,2k}\oplus z_{2b,2k}=Cx_{2k}+v_{2k}+q_{2k+1}$,
\begin{align*}
\hat x_{2k|2k+1}&=\hat x_{2k|2k-1}+L_{2k+1}(z_{2k+1}-C\hat x_{2k|2k-1}),
\end{align*}
\begin{align*}
\begin{bmatrix}x_{2k}\\\hat x_{2k|2k}\\\hat x_{2k|2k+1}\end{bmatrix}
&=\begin{bmatrix}0&I&0\\0&0&I\\(I-L_{2k+1}C)&L_{2k+1}C&0\end{bmatrix}
\begin{bmatrix}\hat x_{2k|2k-1}\\x_{2k}\\\hat x_{2k|2k}\end{bmatrix}\\
&+\begin{bmatrix}0\\0\\L_{2k+1}\end{bmatrix}v_{2k}
+\begin{bmatrix}0\\0\\L_{2k+1}\end{bmatrix}q_{2k+1}.
\end{align*}
State and state estimate update.
\begin{align*}
x_{2k+1}&=Ax_{2k}-BK\hat x_{2k|2k}+w_{2k},\\
\hat x_{2k+1|2k+1}&=A\hat x_{2k|2k+1}-BK\hat x_{2k|2k}.
\end{align*}
denote,
\begin{align*}
\begin{bmatrix}x_{2k}\\\hat x_{2k|2k}\\x_{2k+1}\\\hat x_{2k+1|2k+1}\end{bmatrix}
&=\begin{bmatrix}I&0&0\\0&I&0\\A&-BK&0\\0&-BK&A\end{bmatrix}
\begin{bmatrix}x_{2k}\\\hat x_{2k|2k}\\\hat x_{2k|2k+1}\end{bmatrix}\nonumber\\
&\hskip 10 mm+\begin{bmatrix}0\\0\\I\\0\end{bmatrix}w_{2k}.
\end{align*}
Now define
\begin{align*}
\lim_{k\to \infty}L_{2k}=L_{\text{even}}=\Sigma^{p,\infty}_{II} C^T\left(C\Sigma^{p,\infty}_{II} C^T+R+S_b\right)^{-1},\\
\lim_{k\to \infty}L_{2k+1}=L_{\text{odd}}=\Sigma^{p,\infty}_{II}C^T(C\Sigma^{p,\infty}_{II}C^T+R+S_{2b})^{-1},
\end{align*}
\begin{align*}
\begin{array}{ll}
F_1=\begin{bmatrix}0&0&0&A-BK\\0&0&A&-BK\end{bmatrix},&\hskip -9mm
F_2=\begin{bmatrix}I&0\\0&I\\(I-L_{\text{even}}C)&L_{\text{even}}C\end{bmatrix},\\\\
F_3=\begin{bmatrix}0&I&0\\0&0&I\\(I-L_{\text{odd}}C)&L_{\text{odd}}C&0\end{bmatrix},&\hskip -2mm
F_4=\begin{bmatrix}I&0&0\\0&I&0\\A&-BK&0\\0&-BK&A\end{bmatrix},\\\\
\end{array}\nonumber
\end{align*}
\begin{align}
\begin{array}{llll}
G_1=\begin{bmatrix}0\\I\end{bmatrix},&
G_2=\begin{bmatrix}0\\0\\L_{\text{even}}\end{bmatrix},&
G_3=\begin{bmatrix}0\\0\\L_{\text{odd}}\end{bmatrix},&
G_4=\begin{bmatrix}0\\0\\I\\0\end{bmatrix}
\end{array}\nonumber
\end{align}
Then, the two-step update is described by the recursion
\begin{align*}
\begin{bmatrix}x_{2k}\\\hat x_{2k|2k}\\x_{2k+1}\\\hat x_{2k+1|2k+1}\end{bmatrix}
&=\mathcal M_2\begin{bmatrix}x_{2k-2}\\\hat x_{2k-2|2k-2}\\x_{2k-1}\\\hat x_{2k-1|2k-1}\end{bmatrix}
+\mathcal N_2\begin{bmatrix}w_{2k-1}\\w_{2k}\\\begin{bmatrix}v_{2k}+q_{2k}\\v_{2k}+q_{2k+1}\end{bmatrix}
\end{bmatrix},
\end{align*}
with
\begin{align*}
\mathcal M_2=F_4F_3F_2F_1,\;\;\;
\mathcal N_2=\begin{bmatrix}F_4F_3F_2G_1&G_4&F_4F_3G_2&F_4G_3\end{bmatrix}.
\end{align*}
Whence,
\begin{small}
\begin{align*}
\Psi_{II}=E\left(\begin{bmatrix}x_{2k}\\\hat x_{2k|2k}\\x_{2k+1}\\\hat x_{2k+1|2k+1}\end{bmatrix}
\begin{bmatrix}x^T_{2k}&\hat x^T_{2k|2k}&x^T_{2k+1}&\hat x^T_{2k+1|2k+1}\end{bmatrix}\right),
\end{align*}
\begin{align*}
\Psi_{II}=\texttt{dlyap}\left(\mathcal M_2,\mathcal N_2
\begin{bmatrix}Q&0&0&0\\0&Q&0&0\\0&0&R+S_{b}&R+S_{2b}\\0&0&R+S_{2b}&R+S_{2b}\end{bmatrix}\mathcal N^T_1\right),
\end{align*}
\end{small}
and
\begin{align*}
J_{II}&=\frac{1}{2}\text{trace}\left\{Q_c[\Psi_{II}(1,1)+\Psi_{II}(3,3)]\right\}+\nonumber\\     
    &\hskip 10mm  \frac{1}{2}\text{trace}\left\{K^TR_cK[\Psi_{II}(2,2)+\Psi_{II}(4,4)]\right\},
\end{align*}
%
\emph{Proof for Theorem \ref{the:the3}:}\\
From $\Sigma=\Sigma_{2k|2k-1}$, compute the filter gains
\begin{align*}
L_{2k}&=\Sigma_{2k|2k-1} C^T(C\Sigma_{2k|2k-1} C^T+R+S_b)^{-1},\\
L^\prime_{2k}&=\Sigma_{2k|2k-1} C^T(C\Sigma_{2k|2k-1} C^T+R+S_{b+r})^{-1},\\
\Sigma_{2k+1|2k}&=A\Sigma^\prime_{2k|2k} A^T+Q,\nonumber\\
&=A\Sigma A^T-A\Sigma C^T(C\Sigma C^T+R+S_{b+r})^{-1}C\Sigma A^T+Q,\nonumber\\
L_{2k+1}&=\Sigma_{2k+1|2k}C^T(C\Sigma_{2k+1|2k}C^T+S_{b-r})^{-1}.
\end{align*}
%
 and then,
\begin{align*}
x_{2k}&=Ax_{2k-1}-BK\hat x_{2k-1|2k-1}+w_{2k-1},\\
\hat x_{2k|2k-1}&=A\hat x_{2k-1|2k-1}-BK\hat x_{2k-1|2k-1},\\
&=(A-BK)\hat x_{2k-1|2k-1},\\
\hat x^1_{2k|2k}&=(I-L_{2k}C)\hat x_{2k|2k-1}+L_{2k}z^1_{2k},\\
&=(I-L_{2k}C)\hat x_{2k|2k-1}+L_{2k}Cx_{2k}\\
&+L_{2k}v_{2k}+L_{2k}q_{b,2k},\\
\hat x^\prime_{2k|2k}&=(I-L^\prime_{2k}C)\hat x_{2k|2k-1}+L^\prime_{2k}Cx_{2k}\\
&+L^\prime_{2k}v_{2k}+L^\prime_{2k}q_{b+r,2k},\\
x_{2k+1}&=Ax_{2k}-BK\hat x^1_{2k|2k}+w_{2k},\\
\hat x_{2k+1|2k}&=A\hat x^\prime _{2k|2k}-BK\hat x^1_{2k|2k},\\
\hat x_{2k+1|2k+1}&=(I-L_{2k+1}C)\hat x_{2k+1|2k}+L_{2k+1}z_{2k+1},\\
&=(I-L_{2k+1}C)A\hat x^\prime_{2k|2k}\nonumber\\
&\hskip 5mm-(I-L_{2k+1}C)BK\hat x^1_{2k|2k}\\
&\hskip 7mm+L_{2k+1}Cx_{2k+1}+L_{2k+1}v_{2k+1}\\
&\hskip 10mm+L_{2k+1}q_{b-r,2k+1}.
\end{align*}
The short sequence.
\begin{small}
\begin{align*}
\begin{bmatrix}x_{2k}\\\hat x_{2k|2k-1}\end{bmatrix}
&=\begin{bmatrix}A&-BK&0&0\\0&A-BK&0&0\end{bmatrix}
\begin{bmatrix}x_{2k-1}\\\hat x_{2k-1|2k-1}\\x_{2k-2}\\\hat x^1_{2k-2|2k-2}\end{bmatrix}\\
&\hskip 10mm+\begin{bmatrix}I\\0\end{bmatrix}w_{2k-1},\\
\begin{bmatrix}x_{2k}\\\hat x^1_{2k|2k}\\\hat x^\prime_{2k|2k}\end{bmatrix}
&=\begin{bmatrix}I&0\\L_{2k}C&(I-L_{2k}C)\\L^\prime_{2k}C&(I-L^\prime_{2k}C)\end{bmatrix}
\begin{bmatrix}x_{2k}\\\hat x_{2k|2k-1}\end{bmatrix}\\%
&+\begin{bmatrix}0&0\\L_{2k}&0\\0&L^\prime_{2k}\end{bmatrix}\begin{bmatrix}v_{2k}+q_{b,2k}\\v_{2k}+q_{b+r,2k}\end{bmatrix},\\
\begin{bmatrix}x_{2k+1}\\\hat x^\prime_{2k|2k}\\x_{2k}\\\hat x^1_{2k|2k}\end{bmatrix}
&=\begin{bmatrix}A&-BK&0\\0&0&I\\I&0&0\\0&I&0\end{bmatrix}
\begin{bmatrix}x_{2k}\\\hat x^1_{2k|2k}\\\hat x^\prime_{2k|2k}\end{bmatrix}
+\begin{bmatrix}I\\0\\0\\0\end{bmatrix}w_{2k},\\
\begin{bmatrix}x_{2k+1}\\\hat x_{2k+1|2k+1}\\x_{2k}\\\hat x^1_{2k|2k}\end{bmatrix}
&=\\
&\hskip -22mm\begin{bmatrix}I&0&0&0\\L_{2k+1}C&(I-L_{2k+1}C)A&0&-(I-L_{2k+1}C)BK\\
0&0&I&0\\0&0&0&I\end{bmatrix}%
\times\begin{bmatrix}x_{2k+1}\\\hat x^\prime_{2k|2k}\\x_{2k}\\\hat x^1_{2k|2k}\end{bmatrix}\\
&\hskip 10mm+\begin{bmatrix}0\\L_{2k+1}\\0\\0\end{bmatrix}(v_{2k+1}+q_{b-r,2k+1}),
\end{align*}
\end{small}
denote
\begin{align*}
\lim_{k\to \infty}L_{2k}=L_{\text{even}}=\Sigma^{p,\infty} C^T\left(C\Sigma^{p,\infty} C^T+R+S_b\right)^{-1},\\
\lim_{k\to \infty}L_{2k}'=L_{\text{odd1}}=\Sigma^{p,\infty}C^T(C\Sigma^{p,\infty}C^T+R+S_{b+r})^{-1},\\
\lim_{k\to \infty}L_{2k+1}=L_{\text{odd2}}=\Sigma^{p,\infty}C^T(C\Sigma^{p,\infty}C^T+R+S_{b-r})^{-1},
\end{align*}
\begin{align*}
F_1=\begin{bmatrix}A&-BK&0&0\\0&A-BK&0&0\end{bmatrix},\hskip 2mm &G_1=\begin{bmatrix}I\\0\end{bmatrix},
\end{align*}
\begin{align*}
F_2=\begin{bmatrix}I&0\\L_{\text{even}}C&(I-L_{\text{even}}C)\\L_{\text{odd1}}C&(I-L_{\text{odd1}}C)\end{bmatrix},
G_2=\begin{bmatrix}0&0\\L_{\text{even}}&0\\0&L_{\text{odd1}}\end{bmatrix},
\end{align*}
\begin{align*}
F_3=\begin{bmatrix}A&-BK&0\\0&0&I\\I&0&0\\0&I&0\end{bmatrix},\hskip 2mm G_3=\begin{bmatrix}I\\0\\0\\0\end{bmatrix},\hskip 2mm G_4=\begin{bmatrix}0\\L_{\text{odd2}}\\0\\0\end{bmatrix},
\end{align*}
\begin{align*}
F_4=\begin{bmatrix}I&0&0&0\\L_{\text{odd2}}C&(I-L_{\text{odd2}}C)A&0&-(I-L_{\text{odd2}}C)BK\\
0&0&I&0\\0&0&0&I\end{bmatrix}.
\end{align*}
Then, the two-step update is described by
\begin{align*}
\begin{bmatrix}x_{2k+1}\\\hat x_{2k+1|2k+1}\\x_{2k}\\\hat x^1_{2k|2k}\end{bmatrix}
&=\mathcal M_3\begin{bmatrix}x_{2k-1}\\\hat x_{2k-1|2k-1}\\x_{2k-1}\\\hat x^1_{2k-2|2k-2}\end{bmatrix}
+\mathcal N_3\begin{bmatrix}w_{2k-1}\\w_{2k}\\\begin{bmatrix}v_{2k}+q_{b,2k}\\v_{2k}+q_{b+r,2k}\end{bmatrix}\\v_{2k+1}+q_{b-r,2k+1}\end{bmatrix},
\end{align*}
with
\begin{align*}
\mathcal M_3&=F_4F_3F_2F_1,\\
\mathcal N_3&=\begin{bmatrix}F_4F_3F_2G_1&F_4G_3&F_4F_3G_2&G_4\end{bmatrix},
\end{align*}
\begin{small}
\begin{align*}
\Psi_{III}=E\left(\begin{bmatrix}x_{2k+1}\\\hat x_{2k+1|2k+1}\\x_{2k}\\\hat x^1_{2k|2k}\end{bmatrix}
\begin{bmatrix}x^T_{2k+1}&\hat x^T_{2k+1|2k+1}&x^T_{2k}&\hat x^{1^T}_{2k|2k}\end{bmatrix}\right),
\end{align*}
\end{small}
\begin{align*}
\Psi_{III}=\texttt{dlyap}\left(\mathcal M_3,\mathcal N_3
\mathcal P_3\mathcal N^T_3\right),
\end{align*}
\begin{align*}
\mathcal P_3=\begin{bmatrix}Q&0&0&0&0\\0&Q&0&0&0\\0&0&R+S_b&R+S_{b+r}&0\\0&0&R+S_{b+r}&R+S_{b+r}&0\\0&0&0&0&R+S_{b-r}\end{bmatrix},
\end{align*}
%
\begin{align*}
J_{III}&=\frac{1}{2}\text{trace}\left\{Q_c[\Psi_{III}(1,1)+\Psi_{III}(3,3)]\right\}+\nonumber\\     
    &\hskip 10mm  \frac{1}{2}\text{trace}\left\{K^TR_cK[\Psi_{III}(2,2)+\Psi_{III}(4,4)]\right\}.
\end{align*}

\emph{ Proof for Theorem \ref{the:esc}:}\\
Suppose Pr$(|z_k|>\zeta)=\beta,$ then Pr$(|z_k|<\zeta)=1-\beta,$ for $t=1,2,...$ and assuming the events to be independent,
\begin{align*}
\text{Pr}[(|z_1|<\zeta)\cap(|z_2|<\zeta)...\cap(|z_N|<\zeta)]=(1-\beta)^N.
\end{align*}
The probability that the process escapes at time $T$ is computed as $(1-\beta)^{T-1}\beta$ and the expected time is
\begin{align*}
\text{E}[\tau_{esc}]&=\beta+2(1-\beta)\beta+3(1-\beta)^2\beta+4(1-\beta)^3\beta+...,\\
&=\beta[1+2(1-\beta)+3(1-\beta)^2+4(1-\beta)^3+...],\\
&=\beta\frac{d}{d\beta}[\frac{-1}{1-(1-\beta)}],\\
&=\beta\frac{d}{d\beta}\left[\frac{-1}{\beta}\right],\\
&=\frac{1}{\beta}.
\end{align*}

\bibliographystyle{IEEEtran}
\bibliography{IEEEabrv,/Users/behrooz/PhDpapers/bob}
\end{document}